\documentclass[12pt]{article}

\usepackage{authblk}
\usepackage{tikz}
\usepackage{subfigure}
\usepackage{pgfplots}
\usetikzlibrary{decorations.markings}
\usetikzlibrary{plotmarks}
\usepackage{graphics}
\usepackage{float}
\usepackage{amssymb}
\usepackage{amsfonts}
\usepackage{amsmath}
\usepackage{amsthm}
\usepackage{pdfsync}
\newtheorem{theorem}{Theorem}
\newtheorem{lemma}{Lemma}
\newtheorem{definition}{Definition}

\newtheorem{example}{Example}

\usepackage{mathtools}
\DeclarePairedDelimiterX\Set[2]{\lbrace}{\rbrace}%
 { #1 \,\delimsize|\, #2 }

\usepackage[normalem]{ulem} 
\pgfplotsset{tick label style={font=\Huge},    label style={font=\Huge},    legend style={font=\Huge}}

\usepackage{amssymb}

\usepackage{amsthm}

\newcommand{\ie}{{\it i.e.},\ }

\title{The Invariant Measure of Random Walks in the Quarter-plane: Representation in Geometric Terms}

\author[1]{Yanting Chen}
\author[1]{Richard J. Boucherie}
\author[1,2]{Jasper Goseling}
\affil[1]{Stochastic Operations Research, University of Twente, The Netherlands}
\affil[2]{Department of Intelligent Systems, Delft University of Technology, The Netherlands}
\affil[ ]{\small\texttt{\{y.chen,r.j.boucherie,j.goseling\}@utwente.nl}}

\begin{document}
\maketitle

\begin{abstract}
We consider the invariant measure of homogeneous random walks in the quarter-plane. In particular, we consider measures that can be expressed as a finite linear combination of geometric terms and present conditions on the structure of these linear combinations such that the resulting measure may yield an invariant measure of a random walk. We demonstrate that each geometric term must individually satisfy the balance equations in the interior of the state space and further show that the geometric terms in an invariant measure must have a pairwise-coupled structure. Finally, we show that at least one of the coefficients in the linear combination must be negative.
\end{abstract}


\section{Introduction} \label{sec:introduction}
We study random walks in the quarter-plane that are homogeneous in the sense that transition probabilities are translation invariant. Our interest is in invariant measures that can be expressed as a linear combination of geometric terms, \ie the measure $m$ in state $(i,j)$ is of the form
\begin{equation}{\label{eq:invariantmeasure}}
 m(i,j) = \sum_{(\rho,\sigma) \in \Gamma} \alpha(\rho, \sigma)\rho^i\sigma^j.
\end{equation}

Random walks for which the invariant measure is a geometric product-form are often used to model practical systems. For example, Jackson networks are used to study real systems, see, e.g.,~\cite[Chapter 6]{wolff1989stochastic}. The benefit of such models is that their performance can be readily evaluated with tractable closed-form expressions. 
The performance of systems that do not have a product-form invariant measure can often be approximated by perturbing the transition probabilities to obtain an product-form invariant measure, see e.g.,~\cite[Chapter 9]{van2010queueing}. Various approaches to obtaining comparison results as well as bounds on the perturbation errors exist in the literature, see,~\cite{van1988perturbation, goseling2012linear, muller2002comparison}.

Even though random walks that have a product-form invariant measure have been successfully used for performance evaluation, this class of random walks is rather restricted~\cite[Chapters 1, 5, 6]{van2010queueing}. As a consequence, in many applications it is often not possible to obtain exact results. Therefore, it is of interest to find larger classes of random walks with a tractable invariant measure. Such classes cannot only be of interest for exact performance analysis, but may also be the bases for improved approximation schemes.

For some random walks the invariant measure can be expressed as a linear combination of countably many geometric terms~\cite{adan1993compensation}. This naturally leads to the problem: What are the properties of invariant measures of random walks that are a linear combination of geometric measures? In this paper, we restrict our attention to measures that are a \emph{linear combination of a finite number of geometric measures}. We present conditions on the structure of these linear combinations such that the resulting measure can be an invariant measure of a random walk in the quarter-plane. Our \emph{contributions} are as follows.

For geometric terms $\rho^i \sigma^j$ contained in the summation in~\eqref{eq:invariantmeasure} such that both $\rho > 0$ and $\sigma > 0$, we obtain the following results: First, we demonstrate that each geometric term must individually satisfy the balance equations in the interior of the state space. Second, it is shown that the geometric terms in an invariant measure must have a pairwise-coupled structure stating that for each $(\rho, \sigma)$ in the summation in~\eqref{eq:invariantmeasure} there exists a $(\tilde{\rho}, \tilde{\sigma})$ such that $\tilde{\rho} = \rho$ or $\tilde{\sigma} = \sigma$. Finally, it is shown that if a finite linear combination of geometric terms is an invariant measure, then at least one coefficient $\alpha(\rho, \sigma)$ in~\eqref{eq:invariantmeasure} must be negative.

Various approaches to finding the invariant measure of a random walk in the quarter-plane exist. Most notably, methods from complex analysis have been used to obtain the generating function of the invariant measure~\cite{cohen1983boundary,fayolle1999random}. Matrix-geometric methods provide an algorithmic approach to finding the invariant measure~\cite{neuts1981matrix}. However, explicit closed form expressions for the invariant measures of random walks are hard to obtain using these methods. An overview of the recent work on the tail analysis of the invariant measure of random walks in the quarter-plane is given in~\cite{miyazawa2011light}. 

For reflected Brownian motion with constraints on the boundary transition probabilities, results similar to those reported in the current paper, are presented in~\cite{dieker2009reflected}, where it is shown that for the invariant measure to be a linear combination of exponential measures, there must be an odd number of terms that are generated by a mating procedure, obtaining a structure that we call pairwise-coupled. The method used for the analysis of the continuous state space Brownian motion, however, cannot be used for the discrete state space random walk. Thus, although our results resemble those of~\cite{dieker2009reflected}, the proof techniques substantially differ.

The remainder of this paper is structured as follows. In Section~\ref{sec:model} we present the model. Possible candidates of geometric terms which can lead to an invariant measure are identified in Section~\ref{sec:candidate}. Necessary conditions on the structure of the set of geometric terms are given in Section~\ref{sec:structure}. Section~\ref{sec:CandS} gives conditions on the signs of the coefficients in the linear combination of geometric terms. Several examples of random walks with finite sum of geometric terms invariant measure are provided in Section~\ref{sec:examples}. In Section~\ref{sec:discussion} we summarize our results and present an outlook on future work.

\section{Model} \label{sec:model}

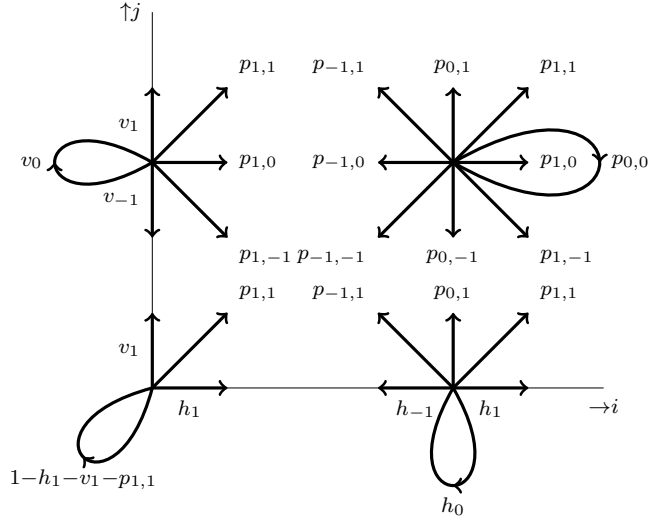
\begin{figure}
\hfill
\begin{tikzpicture}[scale=1]
\tikzstyle{axes}=[very thin] \tikzstyle{trans}=[very thick,->]
\tikzstyle{intloop}=[->,to path={
    .. controls +(30:3) and +(-30:3) .. (\tikztotarget) \tikztonodes}]
\tikzstyle{hloop}=[->,to path={
    .. controls +(-60:2) and +(-120:2) .. (\tikztotarget) \tikztonodes}]
\tikzstyle{vloop}=[->,to path={
    .. controls +(-150:2) and +(-210:2) .. (\tikztotarget) \tikztonodes}]
\tikzstyle{oloop}=[->,to path={
    .. controls +(255:2) and +(195:2) .. (\tikztotarget) \tikztonodes}]

   \draw[axes] (0,0)  -- node[at end, below] {$\scriptstyle \rightarrow i$} (6,0); 
   \draw[axes] (0,0) -- node[at end, left] {$\scriptstyle {\uparrow} {j}$} (0,5);
   \draw[trans] (0,0) to node[below] {$\scriptstyle h_{1}$} (1,0);
   \draw[trans] (0,0) to node [at end, anchor=south west] {$\scriptstyle p_{1,1}$} (1,1);
   \draw[trans] (0,0) to node[left] {$\scriptstyle v_{1}$} (0,1);
   \draw[trans] (4,0) to node[below] {$\scriptstyle h_{-1}$} (3,0);
   \draw[trans] (4,0) to node[below] {$\scriptstyle h_{1}$} (5,0);
   \draw[trans] (4,0) to node[at end, anchor = south west] {$\scriptstyle p_{1,1}$} (5,1);
   \draw[trans] (4,0) to node[at end, anchor = south]  {$\scriptstyle p_{0,1}$} (4,1);
   \draw[trans] (4,0) to node[at end, anchor = south east] {$\scriptstyle p_{-1,1}$} (3,1);
   \draw[trans] (0,3) to node[left] {$\scriptstyle v_{-1}$} (0,2);
   \draw[trans] (0,3) to node[at end, anchor = west] {$\scriptstyle p_{1,0}$} (1,3);
   \draw[trans] (0,3) to node[left] {$\scriptstyle v_{1}$} (0,4);
   \draw[trans] (0,3) to node[at end, anchor = north west] {$\scriptstyle p_{1,-1}$} (1,2);
   \draw[trans] (0,3) to node[at end, anchor = south west] {$\scriptstyle p_{1,1}$} (1,4);
   \draw[trans] (4,3) to node[at end,anchor = west] {$\scriptstyle p_{1,0}$} (5,3);
   \draw[trans] (4,3) to node[at end,anchor = south west] {$\scriptstyle p_{1,1}$} (5,4) ;
   \draw[trans] (4,3) to node[at end, anchor = south] {$\scriptstyle p_{0,1}$} (4,4);
   \draw[trans] (4,3) to node[at end, anchor = south east] {$\scriptstyle p_{-1,1}$} (3,4);
   \draw[trans] (4,3) to node[at end, anchor = east] {$\scriptstyle p_{-1,0}$} (3,3);
   \draw[trans] (4,3) to node[at end, anchor = north east] {$\scriptstyle p_{-1,-1}$} (3,2);
   \draw[trans] (4,3) to node[at end, anchor = north] {$\scriptstyle p_{0,-1}$} (4,2);
   \draw[trans] (4,3) to node[at end, anchor = north west] {$\scriptstyle p_{1,-1}$} (5,2);  
\draw[very thick,decoration={markings, mark=at position 0.5 with {\arrow{>}}},postaction={decorate}] (0,0) to[oloop] node[below] {$\scriptstyle 1 - h_1 - v_1 - p_{1,1} $} (0,0);
\draw[very thick,decoration={markings, mark=at position 0.5 with {\arrow{>}}},postaction={decorate}] (4,0) to[hloop] node[below] {$\scriptstyle h_0$} (4,0);
\draw[very thick,decoration={markings, mark=at position 0.5 with {\arrow{>}}},postaction={decorate}] (4,3) to[intloop] node[right] {$\scriptstyle p_{0,0}$} (4,3);
\draw[very thick,decoration={markings, mark=at position 0.5 with {\arrow{>}}},postaction={decorate}] (0,3) to[vloop] node[left] {$\scriptstyle v_0$} (0,3);
\end{tikzpicture}
\hfill{}
\caption{Random walk in the quarter-plane. \label{fig:rw}}
\end{figure}

Consider a two-dimensional random walk $P$ on the pairs $S = \{(i,j), i,j \in \mathbb{N}_{0}\}$ of non-negative integers. We refer to $\{(i,j) | i>0, j>0\}$, $\{(i,j) | i>0, j=0\}$, $\{(i,j) | i=0, j>0\}$ and $(0,0)$ as the interior, the horizontal axis, the vertical axis and the origin of the state space, respectively. The transition probability from state $(i,j)$ to state $(i+s,j+t)$ is denoted by $p_{s,t}(i,j)$. Transitions are restricted to the adjoining points (horizontally, vertically and diagonally), \ie $p_{s,t}(k,l)=0$ if $|s|>1$ or $|t|>1$. The process is homogeneous in the sense that for each pair $(i,j)$, $(k,l)$ in the interior (respectively on the horizontal axis and on the vertical axis) of the state space it must be that
\begin{equation} \label{eq:homogeneous}
p_{s,t}(i,j)=p_{s,t}(k,l)\quad\text{and}\quad p_{s,t}(i-s,j-t)=p_{s,t}(k-s,l-t),
\end{equation}
for all $-1\leq s\leq 1$ and $-1\leq t\leq 1$.  We introduce, for $i>0$, $j>0$, the notation $p_{s,t}(i,j)=p_{s,t}$, $p_{s,0}(i,0)=h_s$ and $p_{0,t}(0,j)=v_t$.
Note that the first equality of~\eqref{eq:homogeneous} implies that the transition probabilities for each part of the state space are translation invariant. The second equality ensures that also the transition probabilities entering the same part of the state space are translation invariant.
The above definitions imply that $p_{1,0}(0,0)=h_1$ and $p_{0,1}(0,0)=v_1$. 
The model and notation are illustrated in Figure~\ref{fig:rw}. 

We assume that all random walks that we consider are irreducible, aperiodic and positive recurrent. We assume $m$ is the invariant measure, \ie for $i > 0$ and $j > 0$,
\begin{align}
m(i,j) =& \sum_{s = -1}^{1} \sum_{t = -1}^{1} m(i-s,j-t)p_{s,t},\label{eq:interior} \\
m(i,0) =& \sum_{s = -1}^{1} m(i-s,1) p_{s,-1} + \sum_{s = -1}^{1} m(i-s,0)h_s,\label{eq:horizontal} \\
m(0,j) =& \sum_{t = -1}^{1} m(1,j-t) p_{-1,t} + \sum_{t = -1}^{1} m(0,j-t)v_t. \label{eq:vertical}
\end{align}
We will refer to the above equations as the balance equations in the interior, the horizontal axis and the vertical axis, respectively. The balance equation at the origin is implied by the balance equations for all other states.

We are interested in measures that are a linear combination of geometric terms. We first classify the geometric terms.
 
\begin{definition}[{Geometric measures}]{\label{de:degenerateGT}}
The measure $m(i,j) = \rho^i \sigma^j$ is called a geometric measure. It is called horizontally degenerate if $\sigma = 0$, vertically degenerate if $\rho = 0$ and non-degenerate if $\rho > 0$ and $\sigma > 0$. We define $0^0 \equiv 1$.
\end{definition}

We represent a geometric measure $\rho^i \sigma^j$ by its coordinate $(\rho, \sigma)$ in $[0, \infty)^2$. Then, a $\Gamma \subset [0, \infty)^2$ characterizes a set of geometric measures. The set of non-degenerate, horizontally degenerate and vertically degenerate geometric terms from set $\Gamma$ are denoted by $\Gamma_I, \Gamma_H$ and $\Gamma_V$ respectively.

\begin{definition}[{Induced measure}]\label{def:induced}
Signed measure $m$ is called induced by $\Gamma$ if
\begin{equation*}
m(i,j) = \sum_{(\rho, \sigma) \in \Gamma} \alpha(\rho, \sigma)\rho^i\sigma^j,
\end{equation*}
with $\alpha(\rho, \sigma) \in \mathbb{R}\backslash \{0\}$ for all $(\rho, \sigma) \in \Gamma$.
\end{definition}

\begin{figure}
\hfill
\subfigure[\label{fig:2a}]
{
%
%
%
%
\begin{tikzpicture}[scale = 0.2]

\begin{axis}[%
view={0}{90},
width=4.52083333333333in,
height=4.52083333333333in,
scale only axis,
xmin=0, xmax=1.8,
ymin=0, ymax=1.8,
xtick = {-0.5,0,0.5,1,1.5},
ytick = {-0.5,0.5,1,1.5}]

\addplot [draw=red, ultra thick] coordinates{ (0.363636363636364,1.4144017180055)(0.363599457655072,1.41414141414141)(0.358293676367411,1.37373737373737)(0.353552859618717,1.33333333333333)(0.349358154758861,1.29292929292929)(0.345695843347829,1.25252525252525)(0.342557224085144,1.21212121212121)(0.339938618106187,1.17171717171717)(0.337841502046863,1.13131313131313)(0.336272780717225,1.09090909090909)(0.335245218869528,1.05050505050505)(0.334778061114233,1.01010101010101)(0.334897881517189,0.96969696969697)(0.335639721217989,0.929292929292929)(0.337048595630851,0.888888888888889)(0.33918148557815,0.848484848484849)(0.342109973926918,0.808080808080808)(0.345923758612449,0.767676767676768)(0.350735376670035,0.727272727272727)(0.356686632351321,0.686868686868687)(0.363636363636364,0.648281072912731)(0.364069735766821,0.646464646464647)(0.375897026996241,0.606060606060606)(0.390046547588157,0.565656565656566)(0.404040404040404,0.532457127946925)(0.407926181867272,0.525252525252525)(0.434470377019749,0.484848484848485)(0.444444444444444,0.47244094488189)(0.47244094488189,0.444444444444444)(0.484848484848485,0.434470377019749)(0.525252525252525,0.407926181867272)(0.532457127946925,0.404040404040404)(0.565656565656566,0.390046547588157)(0.606060606060606,0.375897026996241)(0.646464646464647,0.364069735766821)(0.648281072912731,0.363636363636364)(0.686868686868687,0.356686632351321)(0.727272727272727,0.350735376670035)(0.767676767676768,0.345923758612449)(0.808080808080808,0.342109973926918)(0.848484848484849,0.33918148557815)(0.888888888888889,0.337048595630851)(0.929292929292929,0.335639721217989)(0.96969696969697,0.334897881517189)(1.01010101010101,0.334778061114233)(1.05050505050505,0.335245218869528)(1.09090909090909,0.336272780717225)(1.13131313131313,0.337841502046863)(1.17171717171717,0.339938618106187)(1.21212121212121,0.342557224085144)(1.25252525252525,0.345695843347829)(1.29292929292929,0.349358154758861)(1.33333333333333,0.353552859618718)(1.37373737373737,0.358293676367411)(1.41414141414141,0.363599457655072)(1.4144017180055,0.363636363636364)(1.45454545454545,0.373575250058153)(1.4949494949495,0.385077253152639)(1.53535353535354,0.398291070708535)(1.55118959054279,0.404040404040404)(1.57575757575758,0.427361969563184)(1.59103880063035,0.444444444444444)(1.58599755101201,0.484848484848485)(1.57575757575758,0.498430462672308)(1.55728453777577,0.525252525252525)(1.53535353535354,0.545805957620431)(1.5158327969253,0.565656565656566)(1.4949494949495,0.582456241811882)(1.46770300443741,0.606060606060606)(1.45454545454545,0.615927537936713)(1.41648436554555,0.646464646464647)(1.41414141414141,0.64817532256267)(1.37373737373737,0.679081653560027)(1.36412783641646,0.686868686868687)(1.33333333333333,0.710211591536339)(1.31204966190001,0.727272727272727)(1.29292929292929,0.741942607491584)(1.26105534769676,0.767676767676768)(1.25252525252525,0.774377885750849)(1.21212121212121,0.807598844810139)(1.21155401326555,0.808080808080808)(1.17171717171717,0.841399360965841)(1.16358095833415,0.848484848484849)(1.13131313131313,0.876430954405704)(1.11743486973948,0.888888888888889)(1.09090909090909,0.912768293900369)(1.07311749366256,0.929292929292929)(1.05050505050505,0.950504218137117)(1.03057927276445,0.96969696969697)(1.01010101010101,0.989740278713788)(0.989740278713788,1.01010101010101)(0.96969696969697,1.03057927276445)(0.950504218137117,1.05050505050505)(0.929292929292929,1.07311749366256)(0.912768293900369,1.09090909090909)(0.888888888888889,1.11743486973948)(0.876430954405704,1.13131313131313)(0.848484848484849,1.16358095833415)(0.841399360965841,1.17171717171717)(0.808080808080808,1.21155401326555)(0.807598844810139,1.21212121212121)(0.774377885750849,1.25252525252525)(0.767676767676768,1.26105534769676)(0.741942607491584,1.29292929292929)(0.727272727272727,1.31204966190001)(0.710211591536339,1.33333333333333)(0.686868686868687,1.36412783641646)(0.679081653560027,1.37373737373737)(0.64817532256267,1.41414141414141)(0.646464646464647,1.41648436554555)(0.615927537936713,1.45454545454545)(0.606060606060606,1.46770300443741)(0.582456241811882,1.4949494949495)(0.565656565656566,1.5158327969253)(0.545805957620431,1.53535353535354)(0.525252525252525,1.55728453777577)(0.498430462672308,1.57575757575758)(0.484848484848485,1.58599755101201)(0.444444444444444,1.59103880063035)(0.427361969563185,1.57575757575758)(0.404040404040404,1.55118959054279)(0.398291070708535,1.53535353535354)(0.385077253152639,1.4949494949495)(0.373575250058153,1.45454545454545)(0.363636363636364,1.4144017180055)(NaN,NaN)};
\addplot [
color=black,
solid,
line width = 1.0pt
]
coordinates{
 (0,1)(1,1) 
};

\addplot [
color=black,
solid,
line width = 1.0pt
]
coordinates{
 (1,0)(1,1) 
};

\end{axis}
\end{tikzpicture}
} 
\subfigure[\label{fig:2b}]
{
%
%
%
%
\begin{tikzpicture}[scale = 0.2]

\begin{axis}[%
view={0}{90},
width=4.52083333333333in,
height=4.52083333333333in,
scale only axis,
xmin=0, xmax=1.4,
ymin=0, ymax=1.4,
xtick = {-0.5,0,0.5,1,1.4},
ytick = {-0.5,0.5,1,1.4}]

\addplot [draw=red, ultra thick] coordinates{ (0.444444444444444,0.855072463768116)(0.4427595386728,0.848484848484849)(0.433506785021937,0.808080808080808)(0.425294122263819,0.767676767676768)(0.418464019350937,0.727272727272727)(0.413527834034268,0.686868686868687)(0.411284562799714,0.646464646464647)(0.413056181651223,0.606060606060606)(0.421199256508875,0.565656565656566)(0.440363228242016,0.525252525252525)(0.444444444444444,0.520430107526882)(0.484848484848485,0.502659574468085)(0.525252525252525,0.503770739064857)(0.565656565656566,0.51387326584177)(0.596857619295741,0.525252525252525)(0.606060606060606,0.529194187582563)(0.646464646464647,0.550306658878505)(0.674585775928315,0.565656565656566)(0.686868686868687,0.573511967064406)(0.727272727272727,0.60031847133758)(0.735808047776655,0.606060606060606)(0.767676767676768,0.631121778387336)(0.787544666332545,0.646464646464647)(0.808080808080808,0.665022421524664)(0.832729455703381,0.686868686868687)(0.848484848484849,0.703258145363409)(0.872043317184383,0.727272727272727)(0.888888888888889,0.747530864197531)(0.905994778322112,0.767676767676768)(0.929292929292929,0.800298661765551)(0.934970443852366,0.808080808080808)(0.958526801448963,0.848484848484849)(0.96969696969697,0.872869318181818)(0.977144955787214,0.888888888888889)(0.99021223901417,0.929292929292929)(0.997881606342452,0.96969696969697)(0.999463517547583,1.01010101010101)(0.994145764871894,1.05050505050505)(0.980964925573209,1.09090909090909)(0.96969696969697,1.11226851851852)(0.957295434623459,1.13131313131313)(0.929292929292929,1.16213674830968)(0.917751097303632,1.17171717171717)(0.888888888888889,1.19086021505376)(0.848484848484849,1.20862649128174)(0.833792470156105,1.21212121212121)(0.808080808080808,1.21738505747126)(0.767676767676768,1.21865430296786)(0.727272727272727,1.21284829721362)(0.725119533896965,1.21212121212121)(0.686868686868687,1.19908998657894)(0.646464646464647,1.17744593881857)(0.639081121252997,1.17171717171717)(0.606060606060606,1.14600113442995)(0.59218447097235,1.13131313131313)(0.565656565656566,1.1031746031746)(0.556985406466246,1.09090909090909)(0.528574445417293,1.05050505050505)(0.525252525252525,1.04569687738005)(0.506784273547258,1.01010101010101)(0.485864447179728,0.96969696969697)(0.484848484848485,0.967647058823529)(0.470608349396229,0.929292929292929)(0.456063907044299,0.888888888888889)(0.444444444444444,0.855072463768116)(NaN,NaN)};
\addplot [
color=black,
solid,
line width = 1.0pt
]
coordinates{
 (0,1)(1,1) 
};

\addplot [
color=black,
solid,
line width = 1.0pt
]
coordinates{
 (1,0)(1,1) 
};

\end{axis}
\end{tikzpicture}
}
\subfigure[\label{fig:2c}] 
{
%
%
%
%
\begin{tikzpicture}[scale = 0.2]

\begin{axis}[%
view={0}{90},
width=4.52083333333333in,
height=4.52083333333333in,
scale only axis,
xmin=0, xmax=1.4,
ymin=0, ymax=1.4,
xtick = {-0.5,0,0.5,1,1.4},
ytick = {-0.5,0.5,1,1.4}]

\addplot [draw=red, ultra thick] coordinates{ (0.202020202020202,0.353682214701604)(0.198384884176416,0.323232323232323)(0.197077895125914,0.282828282828283)(0.202020202020202,0.246508855441227)(0.203210994526056,0.242424242424242)(0.242424242424242,0.203210994526056)(0.246508855441227,0.202020202020202)(0.282828282828283,0.197077895125914)(0.323232323232323,0.198384884176416)(0.353682214701604,0.202020202020202)(0.363636363636364,0.203895429726532)(0.404040404040404,0.214494599816697)(0.444444444444444,0.225785196573688)(0.484848484848485,0.236997817356168)(0.504404788147808,0.242424242424242)(0.525252525252525,0.250700917059662)(0.565656565656566,0.266180848525009)(0.606060606060606,0.28040464748359)(0.613128657834267,0.282828282828283)(0.646464646464647,0.298405096645817)(0.686868686868687,0.315717377758959)(0.705478529499738,0.323232323232323)(0.727272727272727,0.334879782200734)(0.767676767676768,0.355000115222371)(0.786170830398225,0.363636363636364)(0.808080808080808,0.376913130462884)(0.848484848484849,0.399657695376722)(0.856667825816326,0.404040404040404)(0.888888888888889,0.426152914840109)(0.917591296514448,0.444444444444444)(0.929292929292929,0.453964973403925)(0.969690941847002,0.484848484848485)(0.96969696969697,0.484854384740536)(1.01010101010101,0.522600209342491)(1.01309284479204,0.525252525252525)(1.04839442169988,0.565656565656566)(1.05050505050505,0.568865091941115)(1.07584234204985,0.606060606060606)(1.09090909090909,0.636320777503951)(1.09610259029538,0.646464646464647)(1.10924344886468,0.686868686868687)(1.11585119923513,0.727272727272727)(1.11615563320778,0.767676767676768)(1.11040146355617,0.808080808080808)(1.09884567781367,0.848484848484849)(1.09090909090909,0.867567246365329)(1.08131783270928,0.888888888888889)(1.05778644774189,0.929292929292929)(1.05050505050505,0.939687761385969)(1.02741480688696,0.96969696969697)(1.01010101010101,0.989122644265788)(0.989122644265788,1.01010101010101)(0.96969696969697,1.02741480688696)(0.939687761385969,1.05050505050505)(0.929292929292929,1.05778644774189)(0.888888888888889,1.08131783270928)(0.867567246365329,1.09090909090909)(0.848484848484849,1.09884567781367)(0.808080808080808,1.11040146355617)(0.767676767676768,1.11615563320778)(0.727272727272727,1.11585119923513)(0.686868686868687,1.10924344886468)(0.646464646464647,1.09610259029538)(0.636320777503952,1.09090909090909)(0.606060606060606,1.07584234204985)(0.568865091941115,1.05050505050505)(0.565656565656566,1.04839442169988)(0.525252525252525,1.01309284479204)(0.522600209342491,1.01010101010101)(0.484854384740536,0.96969696969697)(0.484848484848485,0.969690941847002)(0.453964973403925,0.929292929292929)(0.444444444444444,0.917591296514448)(0.426152914840109,0.888888888888889)(0.404040404040404,0.856667825816326)(0.399657695376722,0.848484848484848)(0.376913130462884,0.808080808080808)(0.363636363636364,0.786170830398225)(0.355000115222371,0.767676767676768)(0.334879782200734,0.727272727272727)(0.323232323232323,0.705478529499738)(0.315717377758959,0.686868686868687)(0.298405096645817,0.646464646464647)(0.282828282828283,0.613128657834268)(0.28040464748359,0.606060606060606)(0.266180848525009,0.565656565656566)(0.250700917059662,0.525252525252525)(0.242424242424242,0.504404788147808)(0.236997817356168,0.484848484848485)(0.225785196573688,0.444444444444444)(0.214494599816697,0.404040404040404)(0.203895429726532,0.363636363636364)(0.202020202020202,0.353682214701604)(NaN,NaN)};
\addplot [
color=black,
solid,
line width = 1.0pt
]
coordinates{
 (0,1)(1,1) 
};

\addplot [
color=black,
solid,
line width = 1.0pt
]
coordinates{
 (1,0)(1,1) 
};

\end{axis}
\end{tikzpicture}
}
\subfigure[\label{fig:2d}]
{
\input{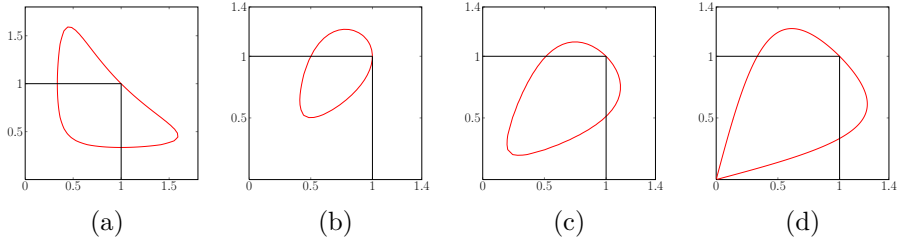}
}
\hfill{}
\caption{Examples of $Q(\rho, \sigma) = 0$. ~\subref{fig:2a} $p_{1,0} = p_{0,1} = \frac{1}{5}$, $p_{-1,-1} = \frac{3}{5}$. ~\subref{fig:2b} $p_{1,0} = \frac{1}{5}$, $p_{0,-1}= p_{-1,1} = \frac{2}{5}$. ~\subref{fig:2c} $p_{1,1} = \frac{1}{62}, p_{-1,1} = p_{1,-1} = \frac{10}{31}, p_{-1,-1} = \frac{21}{62}$. ~\subref{fig:2d} $p_{-1,1} = p_{1,-1} = \frac{1}{4}, p_{-1,-1} = \frac{1}{2}$.}
\label{fig:examplesQ}
\end{figure}

The introduction of signed measures will be convenient in some proofs in Section~\ref{sec:structure}. Our interest is ultimately only in positive measures. If not stated otherwise explicitly, measures are assumed to be positive. To identify the geometric measures that individually satisfy the balance equations in the interior of the state space, \eqref{eq:interior}, we introduce the polynomial
\begin{equation}
 Q(\rho,\sigma) = \rho\sigma\left( \sum_{s = -1}^{1} \sum_{t = -1}^{1} \rho^{-s} \sigma^{-t} p_{s,t} - 1 \right),
\end{equation}
that captures the notion of balance, \ie $Q(\rho,\sigma)=0$ implies that $m(i,j) = \rho^i\sigma^j, i,j \in S$, satisfies~\eqref{eq:interior}.
Several examples of the level sets $Q(\rho, \sigma) = 0$ are displayed in Figure~\ref{fig:examplesQ}.
Let $C$ be the restriction of $Q(\rho,\sigma)=0$ to the interior of the non-negative unit square, \ie
\begin{equation} \label{eq:C}
C = \Big\{ (\rho,\sigma)\in[0,1)^2 \mid
Q(\rho,\sigma)=0
\Big\}.
\end{equation}
In Section~\ref{sec:candidate} we will show that $\Gamma_I \subset C$ is necessary for an induced measure to be the invariant measure of a random walk.

Note that for $|\Gamma|=1$ there are many examples in the literature in which the measure induced by $\Gamma$ is the invariant measure, see, for instance,~\cite[Chapter 6]{wolff1989stochastic}. Also, for $|\Gamma|=\infty$ constructive examples exist, see~\cite{adan1993compensationXXX}. Examples of $\Gamma$ with finite cardinality are provided in Section~\ref{sec:examples}.

\section{Elements in $\Gamma$}{\label{sec:candidate}}

In this section, we obtain conditions on the geometric terms in $\Gamma$ that are necessary for $\Gamma$ to induce an invariant measure of a random walk. We first show that all the non-degenerate geometric terms must come from set $C$. Then we characterize all random walks which may have an invariant measure that includes degenerate geometric terms. Finally we demonstrate that the set $\Gamma$ that induces a measure $m$ is unique.

The next theorem shows that if the measure induced by set $\Gamma$ is the invariant measure, then the non-degenerate geometric terms from set $\Gamma$ must be a subset of $C$, \ie $\Gamma_I \subset C$.

\begin{theorem}{\label{thm:productterm}}
If the invariant measure for a random walk in the quarter-plane is induced by $\Gamma \subset [0,\infty)^2$, where $\Gamma$ is of finite cardinality, then $\Gamma_I \subset C$.
\end{theorem}

We first demonstrate a lemma that will be used in the proof of Theorem~\ref{thm:productterm}.
\begin{lemma}{\label{lem:productterm}}
Let
\[Y = \Set*{n \in \mathbb{N}^+}{\exists (\rho, \sigma) \in \Gamma_I \backslash \{(\rho_1, \sigma_1)\} : \rho_1 \sigma_1^{n} = \rho \sigma^{n}}.\]
Then $|Y| \leq |\Gamma_I| - 1$.
\end{lemma}

\begin{proof}
We will first prove that for any two distinct non-degenerate geometric terms $(\rho_1,\sigma_1)$ and $(\rho,\sigma)$ satisfying $\rho_1 \neq \rho$ and $\sigma_1 \neq \sigma$, there is at most one $n\in \mathbb{N}^+$ for which $\rho_1 \sigma_1^n = \rho \sigma^n$. Assume $\rho_1 \sigma_1^n = \rho \sigma^n$ for some $n \in \mathbb{N}^{+}$. Because $\sigma_1 \neq \sigma$, for any $m \in \mathbb{N}^{+}$ satisfying $m \neq n$, we have $\sigma_1 ^{(m-n)} \neq \sigma^{(m - n)}$. Therefore, $\rho_1 \sigma_1^n \sigma_1^{(m - n)} \neq \rho \sigma^n \sigma^{(m - n)}$, \ie $\rho_1 \sigma_1^m \neq \rho \sigma^m$. From this it follows that there is at most one $n\in \mathbb{N}^+$ for which $\rho_1 \sigma_1^n = \rho \sigma^n$.

It can be readily verified that any non-degenerate geometric term $(\rho, \sigma) \neq (\rho_1, \sigma_1)$ satisfying $\rho = \rho_1$ or $\sigma = \sigma_1$ does not satisfy $\rho_1 \sigma_1^n = \rho \sigma^n$ for any $n \in \mathbb{N}^+$. Moreover, we have shown above that for the non-degenerate geometric term $(\rho, \sigma) \neq (\rho_1, \sigma_1)$ satisfying $\rho \neq \rho_1$ and $\sigma \neq \sigma_1$, there exists at most one positive integer $n$ such that $\rho_1 \sigma_1^{n} = \rho \sigma^{n}$. Therefore, the number of positive integers $n$ for which there exists a $(\rho, \sigma) \in \Gamma_I \backslash \{(\rho_1, \sigma_1)\}$ such that $\rho_1 \sigma_1^{n} = \rho \sigma^{n}$, cannot exceed $|\Gamma_I| - 1$.
\end{proof}

We are now ready to prove Theorem~\ref{thm:productterm}.

\begin{proof}[Proof of Theorem~\ref{thm:productterm}]
Without loss of generality we only prove that $(\rho_1,\sigma_1) \in \Gamma_I$ is in $C$. By deploying Lemma~\ref{lem:productterm}, we conclude that there exists a positive integer $w$ such that for any $(\rho, \sigma) \in \Gamma_I \backslash \{(\rho_1, \sigma_1)\}$, we have $\rho_1 \sigma_1^w \neq \rho \sigma^w$. We now partition $\{(\rho_1, \sigma_1),(\rho_2, \sigma_2), \cdots, (\rho_{|\Gamma_I|}, \sigma_{|\Gamma_I|})\}$ as follows. If $\rho_m \sigma_m^w = \rho_n \sigma_n^w$, then $(\rho_n, \sigma_n)$ and $(\rho_m, \sigma_m)$ will be put into the same element in the partition. We denote this partition by $\Gamma_I^1, \Gamma_I^2, \cdots, \Gamma_I^z$. It is obvious that $(\rho_1, \sigma_1)$ itself form an element and $z \leq |\Gamma_I|$. Without loss of generality, we denote $\Gamma_I^1 = \{(\rho_1, \sigma_1)\}$. Moreover, we arbitrarily choose one geometric term from this element as the representative, which is denoted by $(\rho(\Gamma_I^k), \sigma(\Gamma_I^k))$.

Since the measures induced by $\Gamma_H$ and $\Gamma_V$ are $0$ in the interior of the state space, the balance equation for state $(i,j)$ satisfying $i > 1$ and $j > 1$ is
\begin{equation*}
\sum_{(\rho, \sigma) \in \Gamma_I} \rho^i \sigma^j [\alpha(\rho, \sigma) (1 - \sum_{s = -1}^{1} \sum_{t = -1}^1 \rho^{-s} \sigma^{-t} p_{s,t})] = 0.
\end{equation*}
We now consider the balance equation for states $(d,dw)$ where $d = 2, \cdots, z+ 1$,
\begin{equation*}
\sum_{k = 1}^{z} [\rho(\Gamma_I^k) \sigma(\Gamma_I^k)^w]^d [\sum_{(\rho, \sigma) \in \Gamma_I^k}\alpha(\rho, \sigma) (1 - \sum_{s = -1}^{1} \sum_{t = -1}^1 \rho^{-s} \sigma^{-t} p_{s,t})]= 0.
\end{equation*}
We obtain a system of linear equations in variables $\sum_{(\rho, \sigma) \in \Gamma_I^k}\alpha(\rho, \sigma) (1 - \sum_{s = -1}^{1} \sum_{t = -1}^1 \rho^{-s} \sigma^{-t} p_{s,t})$. The system has a Vandermonde structure in coefficients $[\rho(\Gamma_I^k) \sigma(\Gamma_I^k)^w]^d$. Since any two elements from set
\begin{equation*} 
\{\rho(\Gamma_I^1) \sigma(\Gamma_I^1)^w, \rho(\Gamma_I^2) \sigma(\Gamma_I^2)^w, \cdots, \rho(\Gamma_I^z) \sigma(\Gamma_I^z)^w \}
\end{equation*}
are distinct, we obtain 
\begin{equation*}
1 - \sum_{s = -1}^{1} \sum_{t = -1}^1 \rho_1^{-s} \sigma_1^{-t} p_{s,t} = 0,
\end{equation*}
since $\Gamma_I^1 = \{(\rho_1, \sigma_1)\}$. Therefore, we conclude that $(\rho_1, \sigma_1)$ is in $C$.
\end{proof}

Next, we show that the measure induced by set $\Gamma$ involving degenerate geometric terms cannot be the invariant measure for any random walk.

\begin{theorem}{\label{thm:nice}}
If $\Gamma_H \neq \emptyset$ or $\Gamma_V \neq \emptyset$, then the measure induced by set $\Gamma$ cannot be the invariant measure for any random walk.
\end{theorem}

Before giving the proof of Theorem~\ref{thm:nice}, we provide three technical lemmas. We first give conditions for the sets $\Gamma_H$ and $\Gamma_V$ to be non-empty.
\begin{lemma}{\label{lem:emptygamma}}
If the invariant measure for a random walk in the quarter-plane is
\begin{equation}{\label{eq:completeIM}}
m(i,j) = \sum_{(\rho, \sigma) \in \Gamma_I} \alpha(\rho, \sigma) \rho^i \sigma^j + \sum_{(\rho, 0) \in \Gamma_H} \alpha(\rho, 0) \rho^i 0^j + \sum_{(\sigma, 0) \in \Gamma_V} \alpha(0, \sigma) 0^i \sigma^j,
\end{equation}
then $\Gamma_H \neq \emptyset$ only when $p_{-1,1} = p_{0,1} = p_{1,1} = 0$ and $\Gamma_V \neq \emptyset$ only when $p_{1,-1} = p_{1,0} = p_{1,1} = 0$.
\end{lemma}

\begin{proof}
Since $m(i,j)$ is the invariant measure, $m(i,j)$ satisfies the balance equation at state $(i,1)$ for $i > 1$. Therefore,
\begin{multline}{\label{eq:F}}
\sum_{(\rho, \sigma) \in \Gamma_I} \alpha(\rho, \sigma) \rho^i \sigma =  \sum_{s = -1}^{1} \sum_{t = -1}^{1} \sum_{(\rho, \sigma) \in \Gamma_I} \alpha(\rho, \sigma) \rho^{i-s} \sigma^{1-t} p_{s,t}  +  \\
\sum_{s = -1}^{1} \sum_{(\rho, 0) \in \Gamma_H} \alpha(\rho, 0) \rho^{i-s} p_{s,1}.
\end{multline}
Since $\Gamma_I \subset C$ due to Theorem~\ref{thm:productterm}, equation~\eqref{eq:F} becomes
\begin{equation}{\label{eq:G}}
\sum_{s = -1}^{1} \sum_{(\rho, 0) \in \Gamma_H} \alpha(\rho, 0) \rho^{i-s} p_{s,1} = 0.
\end{equation}
The system of equations for $i = 2, 3, \cdots, |\Gamma_H| + 1$ in equation~\eqref{eq:G} is a Vandermonde system of linear equations if we consider the coefficient $\rho^i$ and unknown $\sum_{s = -1}^{1}  \rho^{-s} p_{s,1}$. Since the elements of $\Gamma_H$ are distinct, we have
\begin{equation}{\label{eq:Hdg}}
\sum_{s = -1}^{1}  \rho^{-s} p_{s,1} = 0.
\end{equation}
for all $(\rho, 0) \in \Gamma_H$. It can be readily verified that only when $\sum_{s = -1}^{1} p_{s,1} = 0$, it is possible to find $\rho \in (0,1)$ such that equation~\eqref{eq:Hdg} is satisfied.
Therefore we conclude that $\Gamma_H$ is non-empty only when $\sum_{s = -1}^{1} p_{s,1} = 0$. Similarly, we conclude that the set $\Gamma_V$ is non-empty only when $\sum_{t = -1}^{1} p_{1,t} = 0$. 
\end{proof}

\begin{lemma}{\label{lem:onlyempty}}
Consider the random walk $P$ in the quarter-plane. {If $m$ induced by set $\Gamma$ is the invariant measure, then $\Gamma_H$ or $\Gamma_V$ must be empty.}
\end{lemma}
\begin{proof}
We know that $\Gamma_H$ is non-empty only when $p_{-1,1} = p_{0,1} = p_{1,1} = 0$ and set $\Gamma_V$ is non-empty only when $p_{1,-1} = p_{1,0} = p_{1,1} = 0$ due to Lemma~\ref{lem:emptygamma}. Assume both $\Gamma_H$ and $\Gamma_V$ are non-empty, we have $p_{-1,1} = p_{0,1} = p_{1,1} = p_{1,0} = p_{1,-1} = 0$, which leads to a reducible random walk. Therefore, we conclude that $\Gamma_H$ or $\Gamma_V$ must be empty.
\end{proof}

The next lemma provides necessary conditions on invariant measure that is induced by $\Gamma$ which includes degenerate geometric terms.
\begin{lemma}{\label{lem:degenerateIM}}
Suppose that the invariant measure for a random walk in the quarter-plane is
\begin{equation}{\label{eq:HIM}}
m(i,j) = \sum_{(\rho, \sigma) \in \Gamma_I} \alpha(\rho, \sigma) \rho^i \sigma^j + \sum_{(\rho, 0) \in \Gamma_H} \alpha(\rho, 0) \rho^i 0^j,
\end{equation}
where set $\Gamma = \Gamma_I \cup \Gamma_H$ is of finite cardinality. Then $m(i,j) = \alpha \rho^i \sigma^j + \tilde{\alpha} \rho^i 0^j$, \ie $\Gamma_I = \{(\rho, \sigma)\}$ and $\Gamma_H = \{(\rho, 0)\}$. Moreover, such a presentation is unique. The result for the invariant measure induced by set $\Gamma = \Gamma_I \cup \Gamma_V$ holds similarly.
\end{lemma}
\begin{proof}
When $\Gamma_I = \emptyset$, the random walk reduces to one dimensional. Hence, we assume $\Gamma_I \neq \emptyset$ here. Since $m(i,j)$ is the invariant measure, $m(i,j)$ satisfies the balance equation for state $(i,0)$ where $i > 1$,
\begin{equation}{\label{eq:H}}
m(i,0) = \sum_{s = -1}^{1} m(i - s,0)h_s + \sum_{s = -1}^{1} m(i - s, 1) p_{s,-1}.
\end{equation}
We will first prove that the invariant measure can only be of the form
\begin{equation}{\label{eq:matchGT}}
m(i,j) = \sum_{k = 1}^{K} (\alpha_{k} \rho_{k}^i \sigma_{k}^j + \tilde{\alpha}_{k} \rho_{k}^i 0^j).
\end{equation}
Substitution of $m(i,j)$ satisfying~\eqref{eq:HIM} in balance equation~\eqref{eq:H} gives
\begin{multline}{\label{eq:nearhorizontal}}
\sum_{(\rho, \sigma) \in \Gamma_I} \alpha(\rho, \sigma) \rho^i(1 - \sum_{s = -1}^{1} \rho^{- s} h_s - \sum_{s = -1}^{1} \rho^{-s} \sigma p_{s,-1}) + \\
\sum_{(\rho, 0) \in \Gamma_H} \alpha(\rho, 0) \rho^i(1 - \sum_{s = -1}^{1} \rho^{ - s} h_s ) = 0.
\end{multline}
Assume there exists a geometric term $(\tilde{\rho}, 0) \in \Gamma_H$ of which the horizontal coordinate is different from that of any geometric terms from set $\Gamma_I$. We now partition set $\Gamma_I \cup \Gamma_H$ as $\Gamma_1, \Gamma_2, \cdots, \Gamma_z$ such that all the geometric terms with the same horizontal coordinates will be put into one element. The common horizontal coordinate is denoted by $\rho(\Gamma_k)$. Clearly, the geometric term $(\tilde{\rho}, 0)$ itself forms an element. Moreover, notice that the non-degenerate geometric term $(\rho, \sigma)$ must satisfy $\sigma = f(\rho)$, where the function $f$ is defined as
\begin{equation}{\label{eq:sbalance}}
f(x) = \frac{1 - (\sum_{s = -1}^{1} x^{-s} p_{s,0})}{\sum_{s = -1}^{1} x^{-s} p_{s,-1}}.
\end{equation} 
Therefore, there is at most one non-degenerate and horizontal degenerate geometric term in set $\Gamma_k$. We now rewrite equation~\eqref{eq:nearhorizontal} as
\begin{multline}{\label{eq:NHX}}
\sum_{k =1}^{z}  \rho(\Gamma_k)^i \sum_{(\rho, \sigma) \in \Gamma_k}[\alpha(\rho, \sigma) (1 - \sum_{s = -1}^{1} \rho^{- s} h_s - \sum_{s = -1}^{1} \rho^{-s} \sigma p_{s,-1})) I[(\rho, \sigma)\in \Gamma_k]  + \\
\alpha(\rho, 0)(1 - \sum_{s = -1}^{1} \rho^{-s} h_s ) I[(\rho, 0) \in \Gamma_k] ] = 0.
\end{multline}
We obtain a system of equations by letting $i = 2, 3, \cdots, |\Gamma_I \cup \Gamma_H| + 1$. This system has a Vandermonde structure by considering the coefficient $\rho(\Gamma_k)$ and the linear relation within the brackets in equation~\eqref{eq:NHX} as unknowns. Since the elements from $\rho(\Gamma_1), \rho(\Gamma_2), \cdots, \rho(\Gamma_z)$ are distinct and the geometric term $(\tilde{\rho}, 0)$ itself forms an element, we obtain
\begin{equation}{\label{eq:tilderho}}
 1 - \sum_{s = -1}^{1} \tilde{\rho}^{- s} h_s  = 0.
\end{equation} 
Because of equation~\eqref{eq:tilderho},  the balance equation~\eqref{eq:nearhorizontal} reduces to
\begin{align}{\label{eq:lessGT}}
\sum_{(\rho, \sigma) \in \Gamma_I} \alpha(\rho, \sigma) \rho^i(1 - \sum_{s = -1}^{1} \rho^{- s} h_s - \sum_{s = -1}^{1} \rho^{-s} \sigma p_{s,-1}) + \notag \\
\sum_{(\rho, 0) \in \Gamma_H \backslash (\tilde{\rho}, 0)} \alpha(\rho, 0) \rho^i(1 - \sum_{s = -1}^{1} \rho^{ - s} h_s ) = 0.
\end{align}
Notice that equation~\eqref{eq:lessGT} is the balance equation for the measure induced by set $\Gamma_I \cup \Gamma_H \backslash (\tilde{\rho}, 0)$. We denote this new measure by $\tilde{m}$. It can be readily verified that measure $\tilde{m}$ is an invariant measure as well. With the same measure in the interior, $m$ has greater measure than $\tilde{m}$ at the horizontal axis, which leads to a contradiction of the uniqueness of the invariant measure for an irreducible ergodic random walk. Similarly, we will draw a contradiction if there exists a geometric term $(\tilde{\rho}, \tilde{\sigma}) \in \Gamma_I$ of which the horizontal coordinate is different from that of any geometric terms from set $\Gamma_H$. Therefore, we have proven that the invariant measure can only be of the form~\eqref{eq:matchGT}. This means the horizontally degenerate geometric terms and non-degenerate geometric terms can only exist in pairs.

Next we will show that $K = 1$ in equation~\eqref{eq:matchGT}. Assume $K > 1$. Without loss of generality we consider a measure $m(i,j)$ with $K = 2$. Since $\Gamma_H \neq \emptyset$ here, we have $\sum_{s = -1}^{1} p_{s,1} = 0$ due to Lemma~\ref{lem:emptygamma}. Moreover, the non-degenerate geometric term $(\rho, \sigma)$ must satisfy $\sigma = f(\rho)$ defined in~\eqref{eq:sbalance}. We observe several properties of $f(x)$. First, $f(x)$ is a continuous function of $x$ and $f(1) = 1$. Secondly, $f(x) = c$ has at most two solutions for any constant $c$. Thirdly, $f(0) \leq 0$. Hence, we conclude that $f(x) = c$ has at most one solution on interval $x \in (0,1)$ when $c \in (0,1)$. This implies that $\rho_1 \neq \rho_2$ and $\sigma_{1} \neq \sigma_{2}$ in measure $m(i,j)$. Moreover, the vertical balance equation for $m(i,j)$ at state $(0,j)$ where $j > 1$ is,
\begin{equation}{\label{eq:vertical2}}
\sum_{k = 1}^{2} \alpha_{k} \sigma_{k}^j(1 - \sum_{t = -1}^{1}\rho_{k}^{- t} v_t - \sum_{t = -1}^{1} \rho_{k}^{-t} \sigma_k p_{-1,t}) = 0.
\end{equation}

We obtain a system of equations when $j = 2,3$. Consider $\sigma_k^j$ as coefficient and $ \alpha_{k} (1 - \sum_{t = -1}^{1}\rho_{k}^{- t} v_t - \sum_{t = -1}^{1} \rho_{k}^{-t} \sigma_k p_{-1,t})$ as unknown, we have a Vandermonde system and therefore obtain that $1 - \sum_{t = -1}^{1}\rho_{k}^{- t} v_t - \sum_{t = -1}^{1} \rho_{k}^{-t} \sigma p_{-1,t} = 0$ for $k = 1,2$. It can be readily verified that both $\alpha_{1} \rho_{1}^i \sigma_{1}^j + \tilde{\alpha}_1 \rho_1^i 0^j$ and $\alpha_{2} \rho_{2}^i \sigma_{2}^j + \tilde{\alpha}_2 \rho_2^i 0^j$ are the invariant measures. Because the invariant measure is unique up to a constant, we have
\begin{equation*} 
\alpha_1 \rho_1^i \sigma_{1}^j = c\alpha_2 \rho_2^i \sigma_2^j,
\end{equation*}
for $i > 1$ and $j > 1$. We obtain a system of equations when $i = 2$ and $j = 2,3$. Consider $\sigma_1^j$, $\sigma_2^j$ as coefficients and $\rho_1^2 \alpha_{1}$, $c \rho_2^2 \alpha_{2}$ as unknowns, we have a Vandermonde system and therefore obtain that $\alpha_k = 0$ for $k = 1,2$, which contradicts the assumption of non-zero coefficients. This also implies that the geometric terms contributed to the invariant measure are unique.
\end{proof}

We are now able to prove Theorem~\ref{thm:nice}.

\begin{proof}[Proof of Theorem~\ref{thm:nice}]
From Lemma~\ref{lem:onlyempty} we know that we cannot have both $\Gamma_{H} \neq \emptyset$ and $\Gamma_{V} \neq \emptyset$. Without loss of generality, let us assume $\Gamma_H \neq \emptyset$. We know from Lemma~\ref{lem:emptygamma} that $p_{-1,1} = p_{0,1} = p_{1,1} = 0$ must be satisfied for the random walk. Therefore, we must have $v_1 > 0$, otherwise the random walk is not irreducible, which violates our assumptions. Moreover, we know from Lemma~\ref{lem:degenerateIM} that if the invariant measure $m(i,j)$ is a sum of geometric terms, it must be of the form $m(i,j) = \alpha \rho^i \sigma^j + \tilde{\alpha} \rho^i 0^j$. Assume $m(i,j)$ is the invariant measure, because $p_{-1,1} = p_{0,1} = p_{1,1}  = 0$, $\tilde{\alpha} \rho^i 0^j$ where $i \geq 0$ and $j \geq 0$ has no contribution to the interior states. Hence, the measure $m_I(i,j) = \alpha \rho^i \sigma^j$ must satisfy the vertical balance~\eqref{eq:vertical}. We now consider the vertical balance equation at state $(0,1)$. Since $m_I(i,j)$ satisfies the vertical balance equation itself, we must have $m_H(i,j) = \tilde{\alpha} \rho^i 0^j$ satisfying the vertical balance equation as well. It can be readily verified that $v_1$ must be zero if $m_H(i,j)$ satisfies the vertical balance equation at state $(0,1)$ for the random walk with $p_{-1,1} = p_{0,1} = p_{1,1} = 0$, hence, we conclude that if $\Gamma_H \neq \emptyset$, then the measure induced by set $\Gamma$ cannot be the invariant measure for any random walk.
\end{proof}

From now on, we restrict ourselves to the non-degenerate geometric terms, \ie $(\rho, \sigma) \in (0,1)^2$.

The next theorem demonstrates that the representation in $\Gamma$ is unique, in the name that adding, deleting or replacing the non-degenerate geometric terms in set $\Gamma$ cannot lead to the same measure $m$. 

\begin{theorem}[Unique representation]{\label{thm:uniquegamma}}
Let $m$ be induced by $\Gamma$ which contains only non-degenerate geometric terms. The representation is unique in the sense that if $m$ is also induced by $\tilde{\Gamma}$, then $\tilde{\Gamma} = \Gamma$.
\end{theorem}

\begin{proof}
Since both $\Gamma$ or $\tilde{\Gamma}$ will lead to $m$, the following equation must hold for all $i > 0$ and $j > 0$,
\begin{equation}{\label{eq:2alpha}}
\sum_{(\rho, \sigma) \in \Gamma \cap \tilde{\Gamma}} (\alpha(\rho, \sigma)-\tilde{\alpha}(\rho, \sigma)) \rho^i \sigma^j + \sum_{(\rho, \sigma) \in \Gamma \backslash \tilde{\Gamma}} \alpha (\rho, \sigma) \rho^i \sigma^j  - \sum_{(\rho, \sigma) \in \tilde{\Gamma} \backslash \Gamma} \tilde{\alpha}(\rho, \sigma) \rho^i \sigma^j = 0.
\end{equation}
We now prove $\alpha(\rho, \sigma) = 0$ for $(\rho, \sigma) \in \Gamma \backslash \tilde{\Gamma}$, $\tilde{\alpha}(\rho, \sigma) = 0$ for $(\rho, \sigma) \in \tilde{\Gamma} \backslash \Gamma$ and $\tilde{\alpha}(\rho, \sigma) = \alpha(\rho, \sigma)$ for $(\rho, \sigma) \in \Gamma \cap \tilde{\Gamma}$. Without loss of generality, we show $\alpha(\rho_1, \sigma_1) - \tilde{\alpha}(\rho_1, \sigma_1) = 0$ for $(\rho_1, \sigma_1) \in \Gamma \cap \tilde{\Gamma}$. Similar to the proof of Theorem~\ref{thm:productterm}, we find a positive integer $w$ and consider a system of equations. This system has a Vandermonde structure with coefficient $(\rho_k \sigma_k^w)^j$ and unknown $\sum_{(\rho, \sigma) \in \Gamma_k} (\alpha(\rho, \sigma) - \tilde{\alpha}(\rho, \sigma))$.  When $(i,j) = (1,w), (2,2w), \cdots, (|\Gamma \cup \tilde{\Gamma}|, |\Gamma \cup \tilde{\Gamma}|w)$, we have a Vandermonde system and obtain that $\tilde{\alpha}(\rho_1, \sigma_1) = \alpha(\rho_1, \sigma_1)$.
\end{proof}

\section{Structure of $\Gamma$} \label{sec:structure}

In this section we consider the structure of $\Gamma$. The proofs in this and the subsequent sections are based on the notion of an uncoupled partition, which is introduced first.  
\begin{definition}[Uncoupled partition]
A partition $\{\Gamma_1, \Gamma_2, \cdots\}$ of $\Gamma$ is \emph{horizontally uncoupled} if $(\rho, \sigma) \in \Gamma_p$ and
$(\tilde{\rho}, \tilde{\sigma}) \in \Gamma_q$ for $p \neq q$, 
implies that $\tilde{\rho} \neq \rho$, \emph{vertically uncoupled}
if $(\rho, \sigma) \in \Gamma_p$ and $(\tilde{\rho}, \tilde{\sigma}) \in \Gamma_q$ for $p \neq q$,
 implies that $\tilde{\sigma} \neq \sigma$, and \emph{uncoupled} if it is both horizontally and vertically uncoupled.
\end{definition}
Horizontally uncoupled sets are obtained by putting pairs $(\rho, \sigma)$ with the same $\rho$ into the same element of the partition. Vertically coupled sets are obtained by putting pairs $(\rho, \sigma)$ with the same $\sigma$ into the same element.

\begin{figure}
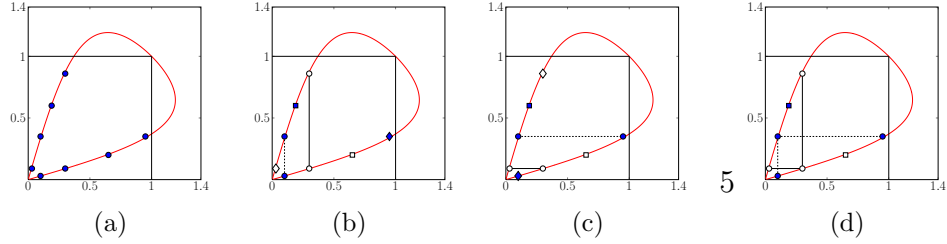

\hfill
\subfigure[
{\label{fig:CandGamma}}
]
{
\input{partition1.tikz}
} 
\subfigure[
 {\label{fig:HUP}}
 ] 
 {
 \input{partition3.tikz}
 }
\subfigure[
{\label{fig:VUP}}
]
  {
  \input{partition4.tikz}
  }
 5 
\subfigure[
{\label{fig:MUP}}
]
{
\input{partition2.tikz}
}
\hfill{}
\caption{Partitions of set $\Gamma$.~\subref{fig:CandGamma} curve $C$ of Figure~\ref{fig:2d} and $\Gamma \subset C$ as dots.~\subref{fig:HUP} horizontally uncoupled partition with $6$ sets.~\subref{fig:VUP} vertically uncoupled partition with $6$ sets.~\subref{fig:MUP} uncoupled partition with $4$ sets. Different sets are marked by different symbols.}
\label{fig:allfigures}
\end{figure} 
We call a partition with the largest number of sets a maximal partition.
\begin{lemma}
The maximal horizontally uncoupled partition, the maximal vertically uncoupled partition and the maximal uncoupled partition are unique.
\end{lemma}
\begin{proof}
Without loss of generality, we only prove that the maximal horizontally uncoupled partition is unique.
Assume that $\{\Gamma_p\}_{p=1}^H$ and $\{\Gamma'_p\}_{p=1}^{H'}$ are different maximal horizontally uncoupled partitions of $\Gamma$. Without loss of generality, $\Gamma_1\cap\Gamma'_1\neq\emptyset$ and $\Gamma_1\setminus\Gamma'_1\neq\emptyset$. Consider $(\rho,\sigma)\in\Gamma_1\setminus\Gamma'_1$ and $(\tilde\rho,\tilde\sigma)\in\Gamma_1\cap\Gamma'_1$. If $\rho=\tilde\rho$, then $\{\Gamma'_p\}_{p = 1}^{H'}$ is not a horizontally uncoupled partition. If $\rho\neq\tilde\rho$, then $\{\Gamma_p\}_{p = 1}^H$ is not maximal. Existence of unique maximal (vertically) uncoupled partitions follows similarly.
\end{proof}
Examples of a maximal horizontally uncoupled partition, of a maximal vertically uncoupled partition and of a maximal uncoupled partition are given in Figure~\ref{fig:allfigures}.
Let $H$ denote the number of elements in the maximal horizontally uncoupled partition and $\Gamma_p^h$, $p=1,\dots,H$, the sets themselves. The common horizontal coordinate of set $\Gamma_p^h$ is denoted by $\rho (\Gamma_p^h)$. The maximal vertically uncoupled partition has $V$ sets, $\Gamma_q^v$, $q=1,\cdots,V$, where elements of $\Gamma_q^v$ have common vertical coordinate $\sigma(\Gamma_q^v)$. The maximal uncoupled partition is denoted by $\{\Gamma_k^u\}_{k=1}^U$.

We start with an observation on the structure of $\Gamma\subset C$ for which the maximal uncoupled partition consists of one set. The degree of $Q(\rho,\sigma)$ is at most two in each variable. Therefore, for each $(\rho,\sigma)\in\Gamma$, there is at most one other geometric term in $\Gamma$ which is horizontally or vertically coupled with $(\rho, \sigma)$. This means, for instance, that if $(\rho,\sigma)\in\Gamma$ and $(\rho,\tilde\sigma)\in \Gamma$, $\tilde\sigma\neq\sigma$, then there does not exist $(\rho,\hat\sigma)\in \Gamma$, where $\hat{\sigma} \neq \sigma$ and $\hat{\sigma} \neq \tilde{\sigma}$. It follows that the elements of $\Gamma$ must be pairwise-coupled.
\begin{definition}[Pairwise-coupled set]
A set $\Gamma\subset C$ is pairwise-coupled if and only if the maximal uncoupled partition of $\Gamma$ contains only one set.
\end{definition}
An example of pairwise-coupled set is 
\begin{equation*}
\Gamma = \{(\rho_k, \sigma_k), k = 1,2,3 \cdots \},
\end{equation*}
 where
\begin{equation*}
 \rho_1 = \rho_2, \sigma_1 > \sigma_2, \rho_2 > \rho_3, \sigma_2 = \sigma_3, \rho_3 = \rho_4, \sigma_3 > \sigma_4, \cdots.
\end{equation*}
%
%
The next theorem states the main result of this section. We show that if there are multiple sets in the maximal uncoupled partition of $\Gamma$,
then the measure induced by this $\Gamma$ cannot be the invariant measure.
\begin{theorem} \label{thm:MUPIM}
Consider the random walk $P$ and its invariant measure $m$. If $m$ is induced by $\Gamma\subset C$, where $\Gamma$ contains only non-degenerate geometric terms, then $\Gamma$ is pairwise-coupled.
\end{theorem}
The proof of the theorem is deferred to the end of this section. We first introduce some additional notation. For any set $\Gamma_p^h$ from the maximal horizontally uncoupled partition of $\Gamma$, let
\begin{align}
B^h(\Gamma_p^h) &= \sum_{(\rho, \sigma)\in\Gamma_p^h}\alpha(\rho, \sigma)[\sum_{s=-1}^1 \big(\rho^{-s} h_s+\rho^{-s}\sigma p_{s,-1}\big) - 1]. \label{eq:Bhdef} 
\end{align}
For any set $\Gamma_q^v$ from the maximal vertically uncoupled partition of $\Gamma$, let
\begin{align}
B^v(\Gamma_q^v) &= \sum_{(\rho, \sigma)\in\Gamma_q^v}\alpha(\rho, \sigma)[\sum_{t=-1}^1 \big(\sigma^{-t} v_t+\rho \sigma^{-t} p_{-1,t}\big) - 1]. \label{eq:Bvdef}
\end{align}
Note that $\sum_{p = 1}^H (\rho(\Gamma_p^h))^{i} B^h(\Gamma_p^h) = 0$ and $\sum_{q = 1}^V (\sigma(\Gamma_q^v))^{j}B^v(\Gamma_q^v) = 0$ are the balance equations for the measure induced by $\Gamma$ at the horizontal and vertical boundary respectively.

The following lemma is a key element for the proof of Theorem~\ref{thm:MUPIM}. It gives the necessary and sufficient conditions for a measure induced by $\Gamma$ to be the invariant measure of a random walk in the quarter-plane.
\begin{lemma} \label{lem:decompbalance}
Consider the random walk $P$ and a measure $m$ induced by $\Gamma\subset C$, where $\Gamma$ contains only non-degenerate geometric terms. Then $m$ is the invariant measure of $P$ if and only if for all $1\leq p\leq H$, $1\leq q\leq V$,
  $B^h(\Gamma^h_p) = 0$ and
  $B^v(\Gamma^v_q) = 0$.
\end{lemma}
\begin{proof}
Since $m$ is the invariant measure of $P$, $m$ satisfies the balance equations at state $(i,0)$. Therefore,
\begin{align} 
0
&= \sum_{s=-1}^1 \big[m(i-s,0)h_s + m(i-s,1)p_{s,-1}\big] - m(i,0) \notag \\
&= \sum_{(\rho, \sigma)\in\Gamma} \alpha(\rho,\sigma) [\sum_{s=-1}^1 \big(\rho^{i-s} h_s+\rho^{i-s}\sigma p_{s,-1}\big) - \rho^i ] \notag  \\
&= \sum_{p=1}^H \rho(\Gamma_p^h)^{i}\sum_{(\rho,\sigma) \in\Gamma_p^h}\alpha(\rho,\sigma)[\sum_{s=-1}^1 \big(\rho^{-s} h_s+\rho^{-s}\sigma p_{s,-1}\big) - 1] \notag \\
&= \sum_{p=1}^H \rho(\Gamma_p^h)^{i} B^h(\Gamma_p^h). \label{eq:A}
\end{align}
From~\eqref{eq:A} it follows that $B^h(\Gamma^h_p)$, $1\leq p\leq H$, satisfy a Vandermonde system of equations. Moreover, from the properties of a maximal horizontally uncoupled partition, the coefficients $\rho(\Gamma^h_p)$ are all distinct. It follows that $B^h(\Gamma^h_p)=0$, $1\leq p\leq H$. Using the same reasoning it follows that $B^v(\Gamma_q^v)=0$, $1\leq q\leq V$, finishing one direction of the proof.

The reversed statement can be verified as follow. If $B^h(\Gamma^h_p)=0$, then $\sum_{p=1}^H (\rho(\Gamma_p^h))^i B^h(\Gamma_p^h) = 0$, where $i = 1,2,3 \cdots$. Therefore, the balance equation for $(i,0)$, $i>0$, is satisfied. Using the same reasoning balance at the vertical states is satisfied. Balance in the interior is satisfied by the assumption that $m$ is induced by $\Gamma\subset C$. Finally, balance in the origin is implied by balance in other parts of the state space.
\end{proof}
We are now ready to present the proof of Theorem~\ref{thm:MUPIM}.
\begin{proof}[Proof of Theorem~\ref{thm:MUPIM}]
The sets of the maximal uncoupled partition can be obtained by taking the union of elements from $\{\Gamma_p^h\}_{p=1}^H$ or $\{\Gamma_q^v\}_{q=1}^V$.
For any $\Gamma^u_k$ where $k=1,\dots,U$, we can find $I_k\subset\{1,\dots,H\}$ and
$J_k\subset\{1,\dots,V\}$ such that $\Gamma_k^u = \bigcup_{p\in I_k} \Gamma_p^h = \bigcup_{q\in J_k} \Gamma_q^v$.
Using the maximal uncoupled partition, we can introduce the signed measures $m_k$, defined as
\begin{equation}
m_k(i,j) = \sum_{(\rho, \sigma)\in\Gamma_k^u} \alpha(\rho, \sigma) \rho^i\sigma^j.
\end{equation}
This allows us to write $m(i,j) = \sum_{k=1}^U m_k(i,j)$. Observe, that $m_k(i,j)$ can be negative.

We will show that if measure $m$ is an invariant measure of the random walk in the quarter-plane, then the measures $m_k$, $k=1,\dots,U,$ will satisfy all balance equations. Let measure $m_k$ be induced by $\Gamma_k$. By the definition of $C$, this implies that all $m_k$, $k=1,\dots,U$, satisfy the balance equations for the states in the interior. Consider the balance equation for $m_k$ at state $(i,0)$. We obtain
\begin{align}
\sum_{s=-1}^1 & [m_k(i-s,0)h_s + m_k(i-s,1)p_{s,-1}] - m_k(i,0) \displaybreak[1] \notag\\
 =&\ \sum_{s=-1}^1 \big[\sum_{\mathclap{(\rho, \sigma)\in\Gamma_k^u}}\alpha(\rho, \sigma)\rho^{i-s} h_s 
   +\sum_{\mathclap{(\rho,\sigma)\in\Gamma_k^u}}\alpha(\rho,\sigma)\rho^{i-s}\sigma p_{s,-1}\big] - \sum_{\mathclap{(\rho, \sigma)\in\Gamma_k^u}}\alpha(\rho, \sigma) \rho^i \displaybreak[1] \notag \\
=&\ \sum_{(\rho, \sigma)\in\Gamma_k^u}\alpha(\rho,\sigma)[\sum_{s=-1}^1 \big(\rho^{i-s} h_s+\rho^{i-s}\sigma p_{s,-1}\big) - \rho^i] \displaybreak[1] \notag \\
=&\ \sum_{p\in I_k}\rho(\Gamma^h_p)^{i}\sum_{(\rho,\sigma)\in\Gamma_p^h}\alpha(\rho,\sigma)[\sum_{s=-1}^1 \big(\rho^{-s} h_s+\rho^{-s}\sigma p_{s,-1}\big) - 1] \displaybreak[1] \notag \\
=&\ \sum_{p\in I_k}\rho(\Gamma^h_p)^{i} B^h(\Gamma_p^h) \notag \displaybreak[1] \\
=&\ 0. \notag
\end{align}
The last equality follows from the assumption that $m$ is an invariant measure and Lemma~\ref{lem:decompbalance}.

In similar fashion it follows that the vertical balance equations of $m_k$ are satisfied as well. As a consequence, we have shown that $m_1, \cdots, m_U$ are signed invariant measures of $P$. Therefore, if $U > 1$ we have a contradiction to Theorem~\ref{thm:uniquegamma} which states the uniqueness of the representation of the sum of geometric terms invariant measure.
\end{proof}

\section{Signs of the coefficients}\label{sec:CandS}
In this section we present conditions on the coefficients $\alpha(\rho, \sigma)$ in the measure induced by $\Gamma$. In particular, we show that at least one of the coefficients in the linear combination must be negative. 
\begin{theorem}\label{thm:negative}
Consider the random walk $P$ and its invariant measure $m$, where $m(i,j) = \sum_{(\rho, \sigma) \in \Gamma} \alpha(\rho, \sigma)\rho^i\sigma^j$, $\Gamma\subset C$, $\alpha(\rho,\sigma) \in \mathbb{R} \backslash \{0\}$. If $m$ is induced by a pairwise-couple set containing only non-degenerate geometric terms, then at least one $\alpha(\rho, \sigma)$ is negative.
\end{theorem}
The proof is based on the following three lemma's. Define
\begin{equation} \label{eq:bhdef}
b^h(\Gamma_p^h) = \frac{B^h(\Gamma_p^h)}{\sum_{(\rho,\sigma)\in \Gamma_p^h }\alpha(\rho,\sigma)} + \left(1-\frac{1}{\rho(\Gamma_p^h)}\right)h_1 + \left(1-\rho(\Gamma_p^h)\right)h_{-1}
\end{equation}
and
\begin{equation} \label{eq:bvdef}
b^v(\Gamma_q^v) = \frac{B^v(\Gamma_q^v)}{\sum_{(\rho,\sigma)\in \Gamma_q^v }\alpha(\rho,\sigma)} + \left(1-\frac{1}{\sigma(\Gamma_q^v)}\right)v_1 + \left(1-\sigma(\Gamma_q^v)\right)v_{-1}.
\end{equation}

\begin{lemma}\label{lem:technicalone}
If $0 < \sigma < \tilde\sigma$, $0 < \rho < \tilde{\rho}$ and {$\alpha(\rho, \sigma) > 0$} then
\begin{align*}
b^h(\{(\rho,\sigma),(\rho,\tilde\sigma)\}) > b^h(\{(\rho,\sigma)\}),\quad & b^h(\{(\rho,\sigma),(\rho,\tilde\sigma)\}) < b^h(\{(\rho,\tilde\sigma)\}), \\
b^v(\{(\rho, \sigma), (\tilde{\rho}, \sigma)\}) > b^v(\{(\rho, \sigma)\}),\quad & b^v(\{(\rho,\sigma),(\tilde\rho,\sigma)\}) < b^v(\{(\tilde\rho,\sigma)\}).
\end{align*}
\end{lemma}
\begin{proof}
From the definition in~\eqref{eq:bhdef} it follows that
\begin{align*}
b^h(\{(\rho, \sigma), (\rho, \tilde{\sigma})\})
= &\frac{\alpha(\rho, \sigma) \sigma + \alpha(\rho, \tilde{\sigma}) \tilde{\sigma}}{\alpha(\rho, \sigma) + \alpha(\rho, \tilde{\sigma})}(\rho p_{-1,-1} + p_{0,-1} + \frac{1}{\rho} p_{1,-1}) - \\
& p_{1,1} - p_{0,1} - p_{-1,1},
\end{align*}
\begin{equation*}
b^h(\{(\rho, \sigma)\}) = \sigma(\rho p_{-1,-1} + p_{0,-1} + \frac{1}{\rho} p_{1,-1}) - p_{1,1} - p_{0,1} - p_{-1,1},
\end{equation*}
and
\begin{equation*}
b^h(\{(\rho, \tilde{\sigma})\}) = \tilde{\sigma}(\rho p_{-1,-1} + p_{0,-1} + \frac{1}{\rho} p_{1,-1}) - p_{1,1} - p_{0,1} - p_{-1,1}.
\end{equation*}
From the above the first row of inequalities follow directly. The remaining inequalities follow directly from~\eqref{eq:bvdef}.
\end{proof}

The following lemma is readily verified and stated without proof.
\begin{lemma}\label{lem:technicaltwo}
If $t_1 (1 - \rho) + t_2 ( 1 - \tilde{\rho}) \geq 0$, $t_1 (1 - 1/\rho) + t_2 (1 - 1/\tilde{\rho}) \geq 0$ and $0 < \rho < \tilde{\rho} < 1$, then $t_1 \leq 0$ and $t_2 \geq 0$.
\end{lemma}

Our final lemma indicates that the linear combination of two non-degenerate geometric terms cannot be the invariant measure of a random walk.
\begin{lemma}\label{lem:finite}
Consider the random walk $P$ and its invariant measure $m$, where $m(i,j) = \sum_{(\rho, \sigma) \in \Gamma} \alpha(\rho, \sigma)\rho^i\sigma^j$, $\Gamma\subset C$, $\alpha(\rho,\sigma) \in \mathbb{R} \backslash \{0\}$. If $m$ is induced by a pairwise-couple set with only non-degenerate geometric terms, then $|\Gamma| \neq 2$.
\end{lemma}
\begin{proof}
Without loss of generality, let
\begin{equation} \label{eq:twopairwisetildem}
m(i,j)=\alpha(\rho, \sigma)\rho^i\sigma^j+\alpha(\rho, \tilde{\sigma})\rho^i\tilde{\sigma}^j,
\end{equation}
where $(\rho, \sigma)\in C$ and $(\rho,\tilde{\sigma})\in C$. It follows from the definition of $C$ that $\sigma$ and $\tilde{\sigma}$ are
the roots of the following quadratic equation in $x$,
\begin{equation}\label{eq:twoquadratic}
\sum_{t = -1}^{1} \sum_{s = -1}^{1} \rho^{-s} p_{s,t} x^{1-t} - x = 0.
\end{equation}

Note that the maximal vertically uncoupled partition of $\{(\rho,\sigma),(\rho,\tilde\sigma)\}$ consists of the two singleton components $\{(\rho,\sigma)\}$ and $\{(\rho,\tilde\sigma)\}$. It follows from Lemma~\ref{lem:decompbalance} that $B^v(\{(\rho,\sigma)\})=B^v(\{(\rho,\tilde\sigma)\})=0$. Therefore, 
$\sigma$ and $\tilde{\sigma}$ are the roots of the following quadratic equation as well
\begin{equation}\label{eq:twobalance}
\sum_{s = -1}^{1} (\rho p_{-1,s} + v_s) x^{1-s} - x = 0.
\end{equation}

From a comparison of the coefficients of~\eqref{eq:twoquadratic} and~\eqref{eq:twobalance} it follows that either \emph{a)} one of the roots will be $1$, contradicting the definition of set $C$ which is restricted within the unit square, or \emph{b)} one geometric term from the pairwise-coupled set must be degenerate. Hence, $m$ cannot be the invariant measure of $P$.
\end{proof}

We are now ready to provide the proof of Theorem~\ref{thm:negative}.
\begin{proof}[Proof of Theorem~\ref{thm:negative}]
Let $(\rho_1, \sigma_1) \in \Gamma$ and $(\rho_2, \sigma_2) \in \Gamma$ satisfy the following conditions:
\begin{itemize}
\item $\rho_1 \geq \rho_2$.
\item $\sigma_1 \geq \sigma_2$.
\item Let $(\rho_1, \sigma_1) \in \Gamma_1^v$, then $\rho_1 \geq \rho$ for all $(\rho, \sigma) \in \Gamma_1^v$.
\item Let $(\rho_1, \sigma_1) \in \Gamma_1^h$, then $\sigma_1 \geq \sigma$ for all $(\rho, \sigma) \in \Gamma_1^h$.
\item Let $(\rho_2, \sigma_2) \in \Gamma_2^v$, then $\rho_2 \leq \rho$ for all $(\rho, \sigma) \in \Gamma_2^v$.
\item Let $(\rho_2, \sigma_2) \in \Gamma_2^h$, then $\sigma_2 \leq \sigma$ for all $(\rho, \sigma) \in \Gamma_2^h$.
\end{itemize}
It can be readily verified that such $(\rho_1, \sigma_1)$, $(\rho_2, \sigma_2)$ always exist.

Without loss of generality, we only discuss the following two cases. In the first case, we have $\rho_1 > \rho_2$ and $\sigma_1 > \sigma_2$. In the second case, we have $\rho_1 = \rho_2$ and $\sigma_1 > \sigma_2$. The proofs for the other cases follow from symmetry considerations.

For the first case we consider the relations
\begin{equation} \label{eq:relcase1}
\begin{aligned}
\left(1 - 1/\rho_1\right)h_1 + (1 - \rho_1)h_{-1} = b^h(\Gamma^h_1), \\
\left(1 - 1/\rho_{2}\right)h_1 + (1 - \rho_{2})h_{-1} = b^h(\Gamma^h_2), \\
\left(1 - 1/\sigma_1\right)v_1 + (1 - \sigma_1)v_{-1} = b^v(\Gamma^v_1), \\
\left(1 - 1/\sigma_{2}\right)v_1 + (1 - \sigma_{2})v_{-1} = b^v(\Gamma^v_2),
\end{aligned}
\end{equation}
which by Lemma~\ref{lem:decompbalance} are required to hold if $m$ is the invariant measure of the random walk $P$. We will construct $s_1$, $s_2$, $t_1$ and $t_2$ that satisfy
\begin{equation} \label{eq:relcase1dual1}
\begin{aligned} 
\left(1 - 1/\rho_1\right)s_1 + \left(1 - 1/\rho_{2}\right)s_2 &\geq 0, \\
\left(1 - \rho_1\right)s_1 + \left(1 - \rho_{2}\right)s_2 &\geq 0, \\
\left(1 - 1/\sigma_1\right)t_1 + \left(1 - 1/\sigma_{2}\right)t_2 &\geq 0, \\
\left(1 - \sigma_1\right)t_1 + \left(1 - \sigma_{2}\right)t_2 &\geq 0 \\
\end{aligned}
\end{equation}
and
\begin{equation} \label{eq:relcase1dual2}
b^h(\Gamma^h_1)s_1 + b^h(\Gamma^h_2)s_2 + b^v(\Gamma^v_1)t_1 + b^v(\Gamma^v_2)t_2 < 0.
\end{equation}
By Farkas' Lemma this leads to a contradiction to~\eqref{eq:relcase1} because the transition probabilities $h_1, h_{-1}, v_1, v_{-1}$ are non-negative. The $s_1$, $s_2$, $t_1$ and $t_2$ are constructed by considering the auxiliary measure $\bar{m} = \alpha(\rho_1,\sigma_1)\rho_1^i \sigma_1^j + \alpha(\rho_{2},\sigma_{2})\rho_{2}^i \sigma_{2}^j$ and the two-dimensional random walk $\bar P$, that has the same transition probabilities as $P$ in the interior of the state space and transition probabilities $\bar h_1$, $\bar h_{-1}$, $\bar v_1$ and $\bar v_{-1}$ along the boundaries. We now consider the relations
\begin{equation} \label{eq:relbar}
\begin{aligned}
\left(1 - 1/\rho_1\right) \bar h_1 + (1 - \rho_1) \bar h_{-1} = b^h(\{(\rho_1,\sigma_1)\}), \\
\left(1 - 1/\rho_{2}\right) \bar h_1 + (1 - \rho_{2}) \bar h_{-1} = b^h(\{(\rho_{2},\sigma_{2})\}), \\
\left(1 - 1/\sigma_1\right) \bar v_1 + (1 - \sigma_1) \bar v_{-1} = b^v(\{(\rho_1,\sigma_1)\}), \\
\left(1 - 1/\sigma_{2}\right) \bar v_1 + (1 - \sigma_{2}) \bar v_{-1} = b^v(\{(\rho_{2},\sigma_{2})\}).
\end{aligned}
\end{equation}
For any non-negative boundary transition probabilities $\bar h_1$, $\bar h_{-1}$, $\bar v_1$ and $\bar v_{-1}$, ~\eqref{eq:relbar} is not satisfied due to Theorem~\ref{thm:MUPIM}. Therefore, by Farkas' Lemma, there exist $s_1$, $s_2$, $t_1$ and $t_2$ that satisfy~\eqref{eq:relcase1dual1} and
\begin{equation*}
b^h(\{(\rho_1, \sigma_1)\})s_1 + b^h(\{(\rho_{2}, \sigma_{2})\})s_2 + b^v(\{(\rho_1, \sigma_1)\})t_1 + b^v(\{(\rho_{2}, \sigma_{2})\})t_2  < 0.
\end{equation*}
Note, that from Lemma~\ref{lem:technicalone} it follows that $b^h(\{\Gamma_1^h\}) \leq b^h(\{(\rho_1, \sigma_1)\})$, $b^h(\{\Gamma_2^h\}) \geq b^h(\{(\rho_{2}, \sigma_{2})\})$,  $b^v(\{\Gamma_1^v\}) \leq b^v(\{(\rho_1, \sigma_1)\})$ and $b^v(\{\Gamma_2^v\}) \geq b^v(\{(\rho_{2}, \sigma_{2})\})$. Also, from Lemma~\ref{lem:technicaltwo} it follows that $s_1\geq 0$, $s_2\leq 0$, $t_1\geq 0$, $t_2\leq 0$. Therefore, $s_1$, $s_2$, $t_1$ and $t_2$ satisfy~\eqref{eq:relcase1dual2}. This concludes the proof of the first case.

For the second case we consider the relations
\begin{equation} \label{eq:relcase2}
\begin{aligned}
\left(1 - 1/\sigma_1\right)v_1 + (1 - \sigma_1)v_{-1} &= b^v(\Gamma^v_1), \\
\left(1 - 1/\sigma_2\right)v_1 + (1 - \sigma_2)v_{-1} &= b^v(\Gamma^v_2),
\end{aligned}
\end{equation}
that are necessary for $m$ to be the invariant measure and obtain a contradiction by constructing $t_1$ and $t_2$ that satisfy
\begin{align} 
\left(1 - 1/\sigma_1\right)t_1 + \left(1 - 1/\sigma_2\right)t_2 &\geq 0, \label{eq:relcase2dual1} \\
\left(1 - \sigma_1\right)t_1 + \left(1 - \sigma_2\right)t_2 &\geq 0, \label{eq:relcase2dual2} \\
b^v(\Gamma^v_1)t_1 + b^v(\Gamma^v_2)t_2 < 0. \label{eq:relcase2dual3}
\end{align}
The auxiliary measure that is used is
$\tilde{m}(i,j) = \alpha(\rho_1,\sigma_1)\rho_1^i \sigma_1^j + \alpha(\rho_2,\sigma_2)\rho_2^i \sigma_2^j$. Observe that $\rho_1=\rho_2$ and that the corresponding relations
are
\begin{align*}
\left(1 - 1/\rho_1\right)h_1 + (1 - \rho_1)h_{-1} &= b^h(\{(\rho_1, \sigma_1), (\rho_2, \sigma_2)\}), \\
\left(1 - 1/\sigma_1\right)v_1 + (1 - \sigma_1)v_{-1} &= b^v(\{(\rho_1, \sigma_1)\}), \\
\left(1 - 1/\sigma_2\right)v_1 + (1 - \sigma_2)v_{-1} &= b^v(\{(\rho_2, \sigma_2)\}).
\end{align*}
From Farkas' Lemma and Lemma~\ref{lem:finite} it follows that there exist $s_1$, $t_1$ and $t_2$ that satisfy~\eqref{eq:relcase2dual1}, \eqref{eq:relcase2dual2} and
\begin{equation}
b^h(\{(\rho_1, \sigma_1), (\rho_2, \sigma_2)\})s_1 + b^v(\{(\rho_1, \sigma_1)\})t_1 + b^v(\{(\rho_2, \sigma_2)\})t_2 \leq 0,
\end{equation}
where $s_1=0$, since it satisfies $(1 - 1/\rho_1)s_1\geq 0$ and $(1 - \rho_1)s_1\geq 0$. Moreover, we have $b^v(\Gamma^v_1)\leq b^v(\{(\rho_1, \sigma_1)\})$ and $b^v(\Gamma^v_2) \geq b^v(\{(\rho_2, \sigma_2)\})$ by Lemma~\ref{lem:technicalone}. In addition, by Lemma~\ref{lem:technicaltwo} we have, $t_1 \geq 0, t_2 \leq 0$. It follows that $t_1$ and $t_2$ satisfy~\eqref{eq:relcase2dual3}.
This concludes the proof of the second case.
\end{proof}

\section{Examples}{\label{sec:examples}}

In this section, we first provide examples of random walks of which the invariant measures are finite mixtures of geometric terms. Then we discuss how such random walks can be constructed. 

The values of the parameters in the examples are mostly obtained as numerical solutions of polynomial equations and are therefore, approximations of the exact results. In addition we depict the transition diagrams of the random walks. In the transition diagrams we have omitted transitions from a state to itself. The examples will be illustrated with a representation of $\Gamma$ on $Q$. In addition of $Q$, we plot in these figures the curves $H$ and $V$ that are the equivalents of $Q$ for the horizontal and vertical balance equations, respectively.

In the first example, we provide a random walk for which the invariant measure is a mixture of three geometric terms. This example also indicates that under favorable conditions, the compensation approach could stop in finitely many steps.  

\begin{example}[Figure~\ref{fig:sq}]{\label{ex:shortestqueue}}
Consider the random walk with $p_{-1,1}=2/5$, $p_{0,-1}=2/5$, $p_{1,-1}=1/5$, $h_1=1/5$, $h_0 = 2/5$, $v_{-1}=18/25$, $v_0 = 2/25$ and all other transition probabilities zero. The measure $m(i,j) = \sum_{k=1}^3 \alpha_k\rho_k^i\sigma_k^j$, where $(\rho_1,\sigma_1)= (1/2, 1/4)$, $(\rho_2,\sigma_2)= (1/16,1/4)$, $(\rho_3,\sigma_3)=(1/16,1/36)$,$\alpha_1 = 1$, $\alpha_2=-20/7$ and $\alpha_3=862/231$ satisfies all balance equations, hence $m(i,j)$ is the invariant measure of the random walk.
\end{example} 

\begin{figure}
\hfill
\subfigure[
{\label{fig:sq1}}
]
{\hfill
\begin{tikzpicture}[scale=1]
\tikzstyle{axes}=[very thin] \tikzstyle{trans}=[very thick,->]
   \draw[axes] (0,0)  -- node[at end, below] {$\scriptstyle \rightarrow i$} (6,0); 
   \draw[axes] (0,0) -- node[at end, left] {$\scriptstyle {\uparrow} {j}$} (0,5);
   \draw[trans] (4,0) to node[below] {$\scriptstyle h_{1}$} (5,0);
   \draw[trans] (4,0) to node[at end, anchor = south east] {$\scriptstyle p_{-1,1}$} (3,1);
   \draw[trans] (0,3) to node[left] {$\scriptstyle v_{-1}$} (0,2);
   \draw[trans] (0,3) to node[at end, anchor = north west] {$\scriptstyle p_{1,-1}$} (1,2);
   \draw[trans] (4,3) to node[at end, anchor = south east] {$\scriptstyle p_{-1,1}$} (3,4);
   \draw[trans] (4,3) to node[at end, anchor = north] {$\scriptstyle p_{0,-1}$} (4,2);
   \draw[trans] (4,3) to node[at end, anchor = north west] {$\scriptstyle p_{1,-1}$} (5,2);  
\end{tikzpicture}
\hfill{}
} 
\subfigure[
 {\label{fig:sq2}}
 ] 
 {
 \input{allcurves_sq.tikz}
 }
\hfill{}
\caption{Example~\ref{ex:shortestqueue}.~\subref{fig:sq1} Transition diagram of Example~\ref{ex:shortestqueue}.~\subref{fig:sq2} Balance equations. The geometric terms contributed to the invariant measure are denoted by the blue squares.\label{fig:sq}}
\end{figure} 

The next example illustrates a random walk with sum of three geometric terms invariant measure without satisfying the constraint $p_{1,0} + p_{1,1} + p_{0,1} = 0$, which is required by compensation approach, see~\cite{adan1993compensation}. This means, for random walks where the compensation approach cannot be applied, the mixture of finite geometric terms invariant measure may still exist.

\begin{example}[Figure~\ref{fig:tdq}]{\label{ex:tandemqueue}}
Consider the random walk with $p_{1,0}=0.05$, $p_{-1,1}=0.15$, $p_{0,-1}=0.15$, $p_{0,0} = 0.65$ $h_1=0.15$, $h_0 = 0.55$, $v_1 = 0.0929$, $v_{-1}=0.15$, $v_0 = 0.7071$ and all other transition probabilities zero. The measure $m(i,j) = \sum_{k=1}^3 \alpha_k\rho_k^i\sigma_k^j$, where $(\rho_1,\sigma_1)= (0.4618,0.3728)$, $(\rho_2,\sigma_2)= (0.2691,0.3728)$, $(\rho_3,\sigma_3)=(0.2691,0.7218)$, {$\alpha_1 = 0.1722$, $\alpha_2=-0.2830$ and $\alpha_3=0.2251$} satisfies all balance equations, hence $m(i,j)$ is the invariant measure of the random walk.
\end{example} 

\begin{figure}
\hfill
\subfigure[
{\label{fig:tdq1}}
]
{\hfill
\begin{tikzpicture}[scale=1]
\tikzstyle{axes}=[very thin] \tikzstyle{trans}=[very thick,->]
   \draw[axes] (0,0)  -- node[at end, below] {$\scriptstyle \rightarrow i$} (6,0); 
   \draw[axes] (0,0) -- node[at end, left] {$\scriptstyle {\uparrow} {j}$} (0,5);
   \draw[trans] (4,0) to node[below] {$\scriptstyle h_{-1}$} (3,0);
   \draw[trans] (4,0) to node[below] {$\scriptstyle h_{1}$} (5,0);
   \draw[trans] (4,0) to node[at end, anchor = south east] {$\scriptstyle p_{-1,1}$} (3,1);
   \draw[trans] (0,3) to node[left] {$\scriptstyle v_{-1}$} (0,2);
   \draw[trans] (0,3) to node[at end, anchor = west] {$\scriptstyle p_{1,0}$} (1,3);
   \draw[trans] (0,3) to node[left] {$\scriptstyle v_{1}$} (0,4);
   \draw[trans] (4,3) to node[at end,anchor = west] {$\scriptstyle p_{1,0}$} (5,3);
   \draw[trans] (4,3) to node[at end, anchor = south east] {$\scriptstyle p_{-1,1}$} (3,4);
   \draw[trans] (4,3) to node[at end, anchor = north] {$\scriptstyle p_{0,-1}$} (4,2);
\end{tikzpicture}
\hfill{}
} 
\subfigure[
 {\label{fig:tdq2}}
 ] 
 {
 \input{allcurves_tdq.tikz}
 }
\hfill{}
\caption{Example~\ref{ex:tandemqueue}.~\subref{fig:tdq1} Transition diagram of Example~\ref{ex:tandemqueue}.~\subref{fig:tdq2} Balance equations. The geometric terms contributed to the invariant measure are denoted by the blue squares.\label{fig:tdq}}

\end{figure} 

The next example uses five geometric terms in the invariant measure.

\begin{example}[Figure~\ref{fig:crbq}]{\label{ex:crbqueue}}
Consider the random walk with $p_{1,0}=0.05$, $p_{0,1} = 0.05$, $p_{-1,1}=0.2$, $p_{-1,0}=0.2$,$p_{0,-1} = 0.2$, $p_{1,-1} = 0.2$, $p_{0,0} = 0.1$, $h_1=0.5$, $h_{-1} = 0.1$, $h_0 = 0.15$, $v_1 = 0.113$, $v_{-1}=0.06$, $v_0 = 0.577$ and all other transition probabilities zero. The measure $m(i,j) = \sum_{k=1}^5 \alpha_k\rho_k^i\sigma_k^j$, where $(\rho_1,\sigma_1)= (0.9773,0.5947)$, $(\rho_2,\sigma_2)= (0.3224,0.5947)$, $(\rho_3,\sigma_3)=(0.3224,0.2346)$, $(\rho_4,\sigma_4)=(0.2857,0.2346)$, $(\rho_5, \sigma_5) = (0.2857,0.5073)$. And $\alpha_1 = 0.0088$, $\alpha_2=0.1180$, $\alpha_3=-0.1557$, $\alpha_4 = 0.1718$, $\alpha_5 = -0.1414$ satisfies all balance equations, hence $m(i,j)$ is the invariant measure of the random walk.
\end{example}

\begin{figure}
\hfill
\subfigure[
{\label{fig:crbq1}}
]
{\hfill
\begin{tikzpicture}[scale=1]
\tikzstyle{axes}=[very thin] \tikzstyle{trans}=[very thick,->]
   \draw[axes] (0,0)  -- node[at end, below] {$\scriptstyle \rightarrow i$} (6,0); 
   \draw[axes] (0,0) -- node[at end, left] {$\scriptstyle {\uparrow} {j}$} (0,5);
   \draw[trans] (4,0) to node[below] {$\scriptstyle h_{-1}$} (3,0);
   \draw[trans] (4,0) to node[below] {$\scriptstyle h_{1}$} (5,0);
   \draw[trans] (4,0) to node[at end, anchor = south]  {$\scriptstyle p_{0,1}$} (4,1);
   \draw[trans] (4,0) to node[at end, anchor = south east] {$\scriptstyle p_{-1,1}$} (3,1);
   \draw[trans] (0,3) to node[left] {$\scriptstyle v_{-1}$} (0,2);
   \draw[trans] (0,3) to node[at end, anchor = west] {$\scriptstyle p_{1,0}$} (1,3);
   \draw[trans] (0,3) to node[left] {$\scriptstyle v_{1}$} (0,4);
   \draw[trans] (0,3) to node[at end, anchor = north west] {$\scriptstyle p_{1,-1}$} (1,2);
   \draw[trans] (4,3) to node[at end,anchor = west] {$\scriptstyle p_{1,0}$} (5,3);
   \draw[trans] (4,3) to node[at end, anchor = south] {$\scriptstyle p_{0,1}$} (4,4);
   \draw[trans] (4,3) to node[at end, anchor = south east] {$\scriptstyle p_{-1,1}$} (3,4);
   \draw[trans] (4,3) to node[at end, anchor = east] {$\scriptstyle p_{-1,0}$} (3,3);
   \draw[trans] (4,3) to node[at end, anchor = north] {$\scriptstyle p_{0,-1}$} (4,2);
   \draw[trans] (4,3) to node[at end, anchor = north west] {$\scriptstyle p_{1,-1}$} (5,2);  
\end{tikzpicture}
\hfill{}
} 
\subfigure[
 {\label{fig:crbq2}}
 ] 
 {
 \input{allcurves_crossb.tikz}
 }
\hfill{}
\caption{Example~\ref{ex:crbqueue}.~\subref{fig:crbq1} Transition diagram of Example~\ref{ex:crbqueue}.~\subref{fig:crbq2} Balance equations. The geometric terms contributed to the invariant measure are denoted by the blue squares.\label{fig:crbq}}

\end{figure}

The construction of a random walk with sum of finite geometric terms invariant measure depends on the locations of the intersections of the boundary balance equations and interior balance equation. If there exists a pairwise-coupled set connecting the intersection of $H$ with $Q$ to the intersection of $V$ with $Q$, then there exists mixture of finite geometric terms invariant measure. We conclude that choosing proper boundary transition probabilities is essential for the existence of sum of finite geometric terms invariant measure.

\section{Conclusion}\label{sec:discussion}
In this paper, we have obtained necessary conditions on measures induced by geometric terms that are the invariant measure of a random walk. In particular, non-degenerate terms must each satisfy the balance equations in the interior of the state space, and must form a pairwise-coupled set. In the linear combination of non-degenerate terms, at least one coefficient must be negative. We have completed the necessary conditions by also including degenerate terms.

It is interesting to note that the pairwise-coupled structure obtained in this paper is equal to the structure obtained in the compensation approach by Adan et al.~\cite{adan1993compensation}. It is suggested in~\cite{adan1993compensation} that the compensation approach, in favorable conditions, might provide a finite number of terms. Our example~\ref{ex:shortestqueue} in Section~\ref{sec:examples} provides a constructive example of such a random walk. Note, however, that the compensation approach, in general, generates countably many geometric terms. It is of interest to generalize the necessary conditions of this paper to the case of countably infinitely many geometric terms. This will require a complete characterization of the algebraic properties of $Q(\rho, \sigma) = 0$ similar to, for instance, the work in~\cite{fayolle1999random}. Since these techniques are fundamentally different from the ones used in the current paper, a generalization to the case of countably infinitely many terms is among our aim for future research.  As part of further work we will also study corresponding sufficient conditions and approximations schemes based on sums of geometric terms. 

Among other possible directions for future research are an extension to higher dimensional walks and random walks with different transition structure, e.g. by allowing longer jumps. The extension of our results to higher dimensional random walks seems feasible using the techniques that we have developed in the current paper. The extension to longer jumps, however, will require substantially different techniques. The reason is that in the current work we have made extensive use of the fact that short jumps induce balance equations that are polynomials of at most degree two.

\section{Acknowledgments}
The authors thank the anonymous reviewer for the useful suggestions. Yanting Chen acknowledges support by a CSC scholarship [No.2008613008]. This work is partly supported by the Netherlands Organisation for Scientific Research (NWO) grant $612.001.107$.


\end{document}